\documentclass[preprint,11pt]{elsarticle}
\usepackage{lmodern}
\usepackage{amsmath,amssymb,amsthm}
\usepackage[margin=1in]{geometry}
\usepackage{microtype}
\usepackage{hyperref}
\usepackage{enumitem}
\usepackage[numbers]{natbib}
\usepackage{amsmath,amssymb,amsthm}
\usepackage{microtype}
\usepackage{hyperref}
\usepackage{enumitem}
\IfFileExists{stix.sty}{}{\usepackage{lmodern}}

\numberwithin{equation}{section}
\newtheorem{theorem}{Theorem}[section]
\newtheorem{lemma}[theorem]{Lemma}
\newtheorem{proposition}[theorem]{Proposition}
\newtheorem{corollary}[theorem]{Corollary}
\theoremstyle{remark}
\newtheorem{remark}[theorem]{Remark}
\newcommand{\Ahat}{\widehat A}
\DeclareMathOperator{\sech}{sech}

\begin{document}
	\begin{frontmatter}
		\title{Log-concavity for alternating Eulerian polynomials}
		
		\author[a]{Lijie Wang}
		\ead{lijiewang26@hotmail.com}
		
		\author[b]{Sainan Zheng\corref{cor1}}
		\ead{zhengsainan@dufe.edu.cn}
		\cortext[cor1]{Corresponding author}
		
		\address[a]{School of Mathematical Sciences, Dalian University of Technology,
			Dalian 116024, P. R. China}
		
		\address[b]{School of Data Science and Artificial Intelligence,
			Dongbei University of Finance and Economics,
			Dalian 116025, P. R. China}
		
		\begin{abstract}
			Alternating Eulerian polynomials enumerate permutations by the
			alternating descent statistic introduced by Chebikin.
			We prove that their coefficient sequences are log-concave,
			verifying a conjecture of Lin, Ma, Wang and Wang.
			Our proof combines a factorization involving reciprocal quadratic
			factors, estimates for ordinary Eulerian zeros, and a log-concavity
			criterion for products of such factors.
			The required finite verifications use exact integer arithmetic.
			As applications, we obtain log-concavity for the distribution of
			$3$-descents, a local limit theorem for alternating descents,
			and bounds for the largest coefficient.
			We also give a Lorentzian reformulation of the main theorem.
		\end{abstract}
		\begin{keyword}
			alternating Eulerian polynomial \sep log-concavity \sep reciprocal quadratic factors \sep Eulerian zeros
			\MSC[2020] 05A20\sep 05A05\sep 26C10
		\end{keyword}
	\end{frontmatter}

\section{Introduction}
Let $\mathfrak S_n$ denote the set of all permutations
of $[n]=\{1,2,\ldots,n\}$.
For a permutation
$\pi=\pi_1\pi_2\cdots\pi_n\in\mathfrak S_n$,
a descent is an index $i\in[n-1]$ such that
$\pi_i>\pi_{i+1}$.
As a variation of the classical descent statistic,
Chebikin~\cite{Chebikin2008} introduced the
alternating descent statistic.
An index $i\in[n-1]$ is an alternating descent
of $\pi$ if $\pi_i>\pi_{i+1}$ when $i$ is odd,
or $\pi_i<\pi_{i+1}$ when $i$ is even.
Let $\operatorname{altdes}(\pi)$ denote the number
of alternating descents of $\pi$.

The \emph{alternating Eulerian polynomial} is defined by
\[
\Ahat_n(t)
=\sum_{\pi\in\mathfrak S_n}
t^{\operatorname{altdes}(\pi)}
=\sum_{k=0}^{n-1}a_{n,k}t^k,
\]
where $a_{n,k}$ counts permutations in
$\mathfrak S_n$ with exactly $k$ alternating
descents, and $\Ahat_n(1)=n!$.
The first few alternating Eulerian polynomials are
$
\Ahat_1(t)=1,
\Ahat_2(t)=1+t,
\Ahat_3(t)=2+2t+2t^2.
$

A polynomial $f(t)=\sum_{k=0}^{n}a_kt^k$
is called \emph{palindromic} if
$a_k=a_{n-k}$ for $0\le k\le n$. 
Replacing each entry $\pi_i$ by $n+1-\pi_i$ reverses every
comparison and changes the number of alternating descents
from $k$ to $n-1-k$.
Consequently, $\Ahat_n(t)$ is palindromic of degree $n-1$;
equivalently,
\begin{equation}\label{eq:symmetry}
a_{n,k}=a_{n,n-1-k}.
\end{equation}

Alternating Eulerian polynomials are related to tangent
derivative polynomials, and their zeros lie on the unit
circle~\cite{Chow2024,MaYeh2016}.
Lin, Ma, Wang and Wang~\cite{LinMaWangWang2022}
studied their sign-alternating $\gamma$-expansions,
and Zhu~\cite{Zhu2018} proved coefficientwise
$2$-log-convexity as the order varies.
Here we study log-concavity of the coefficients
for each fixed order.
Further studies include type $B$ analogues and grammatical
descriptions~\cite{MaFangMansourYeh2022},
combinatorial formulas and parity
refinements~\cite{Pan2022,PanZeng2023},
and $cd$-indices~\cite{PanQian2024}.
Heerten, Sambale and Th\"ale~\cite{HeertenSambaleThale2024}
studied probabilistic limits for alternating descents.

Let $a_0,a_1,\ldots,a_n$ be a sequence of nonnegative real numbers.
The sequence is said to be \emph{unimodal} if there exists an index
$0\le m\le n$ such that
$a_0\le\cdots\le a_{m-1}\le a_m\ge a_{m+1}\ge\cdots\ge a_n.$
The sequence is said to be \emph{log-concave} if
$a_k^2\ge a_{k-1} a_{k+1}$ for $k=1,\ldots,n-1$.
For a polynomial $f(t)=\sum_{k=0}^n a_kt^k$, define
\[
\Delta_k(f)=a_k^2-a_{k-1}a_{k+1},
\]
where coefficients outside the range $0\le k\le n$
are taken to be zero.
Thus, the coefficient sequence of $f(t)$ is
log-concave if and only if
$\Delta_k(f)\ge0$ for $1\le k\le n-1$.
A sequence $(a_k)_{k=0}^n$ has
\emph{no internal zeros} if there are no indices
$0\le i<j<k\le n$ such that
$a_i\ne0$, $a_j=0$, and $a_k\ne0$.
A nonnegative log-concave sequence with no internal zeros is unimodal (see, e.g., Brenti \cite{Brenti1989}).

We use two standard facts.
The convolution property below follows from Brenti~\cite[Theorem~2.2.2 and Proposition~2.5.1]{Brenti1989}.
The mixed coefficient inequality is an adjacent-index case of the inequalities of Ellard and \v{S}migoc~\cite[Theorem~2.2 and Remark~1]{EllardSmigoc2020}.

\begin{lemma}[{\cite[Theorem~2.2.2 and
Proposition~2.5.1]{Brenti1989}}]\label{lem:convolution}
Let $f(t)$ and $g(t)$ be two nonzero polynomials
with nonnegative log-concave coefficient sequences
having no internal zeros.
Then the coefficient sequence of $f(t)g(t)$
is also nonnegative, log-concave, and has
no internal zeros.
\end{lemma}
\begin{lemma}[{
\cite[Theorem~2.2 and Remark~1]{EllardSmigoc2020}}]\label{lem:mixed}
Let
\[
q(t)=\prod_{i=1}^{m}(1+c_it+t^2)
=\sum_k b_kt^k,
\qquad c_i\ge0.
\]
Then
\begin{equation}\label{eq:mixed-products}
b_kb_{k-1}\ge b_{k+1}b_{k-2}
\qquad(k\in\mathbb Z),
\end{equation}
where coefficients outside the range of $q(t)$
are taken to be zero.
\end{lemma}

Since $c_i\ge0$, factoring the quadratics into linear
factors with constant term $1$ gives a self-conjugate
list of parameters in the closed right half-plane.
The coefficients $b_k$ are the corresponding elementary
symmetric functions, so~\eqref{eq:mixed-products} follows
from the adjacent-index inequalities in the cited results.

Lin, Ma, Wang and Wang
\cite[Theorem~1.1 and Conjecture~5.1]{LinMaWangWang2022}
proved the unimodality of the alternating Eulerian
polynomials and conjectured that their coefficient
sequences are log-concave for all $n\ge1$.
We prove this conjecture in the following theorem.

\begin{theorem}\label{thm:main}
For every $n\ge1$, the coefficient sequence of
the alternating Eulerian polynomial $\Ahat_n(t)$
is log-concave.
\end{theorem}

The rest of this paper is organized as follows.
In Section~2, we establish the factorization of
alternating Eulerian polynomials and a lower bound
for the number of factors with parameters exceeding
$3/2$.
In Section~3, we prove an endpoint reduction criterion
and a smoothing lemma.
In Section~4, we establish a block criterion
and complete the proof of Theorem~\ref{thm:main}.
Finally, Section~5 gives applications
to the distribution of $3$-descents and local limits,
bounds for the largest coefficient, and a Lorentzian reformulation of the main theorem.
\section{Factorization and strong-factor estimates} \label{sec:factorization}

The reciprocal-quadratic factorization
follows from the zero results of
Ma and Yeh~\cite{MaYeh2016}.
We identify its parameters in terms
of ordinary Eulerian zeros, which allows
us to estimate the number of strong factors
needed for the block criterion in Section~4.

Let
\[
A_n(z)=\sum_{\pi\in\mathfrak S_n}
z^{\operatorname{des}(\pi)}
\]
be the ordinary Eulerian polynomial, where
$\operatorname{des}(\pi)$ denotes the number of descents.
It is well known that all zeros of $A_n(z)$ are
simple and negative for $n\ge2$; see~\cite{Chow2024}.
The complement map also gives the palindromic identity
\begin{equation}\label{eq:eulerian-pal}
A_n(z)=z^{n-1}A_n(z^{-1}).
\end{equation}
Hence the roots other than $-1$ occur in reciprocal pairs, with one root less than $-1$ and the other in $(-1,0)$.
The root $-1$ occurs precisely when $n$ is even.
We now use these properties to obtain a factorization
of the alternating Eulerian polynomials.

\begin{proposition}\label{prop:factor}
For $n\ge1$, let $m=\lfloor(n-1)/2\rfloor$.
There exist distinct numbers $v_1,\ldots,v_m>0$
such that $A_n(-e^{v_i})=0$ and
\begin{equation}\label{eq:factor}
\Ahat_n(t)
=a_{n,0}(1+t)^{n-1-2m}
\prod_{i=1}^{m}(1+c_it+t^2),
\qquad
c_i=2\sech v_i\in(0,2),
\end{equation}
where $a_{n,0}>0$.
\end{proposition}

\begin{proof}
The cases $n=1,2$ follow from
$\Ahat_1(t)=1$ and $\Ahat_2(t)=1+t$.
Suppose that $n\ge3$.
Let $T_n$ be the tangent derivative polynomial,
defined by
\[
\frac{d^n}{dx^n}\tan x=T_n(\tan x).
\]
Ma and Yeh~\cite{MaYeh2016} established the identity
\begin{equation}\label{eq:mayeh}
2^n(1+t^2)\Ahat_n(t)
=(1-t)^{n+1}
T_n\left(\frac{1+t}{1-t}\right).
\end{equation}

We first relate the roots of $T_n$ to those
of the ordinary Eulerian polynomial $A_n$.
By the classical Eulerian generating-function
identity~\cite[Proposition~1.4.4]{Stanley2012},
after normalization, we have,
\begin{equation}\label{eq:ordinary-derivative}
\left(z\frac{d}{dz}\right)^n
\frac{1}{1-z}
=
\frac{zA_n(z)}{(1-z)^{n+1}}.
\end{equation}

Consider a pair of reciprocal roots
$-\alpha$ and $-\alpha^{-1}$ of $A_n(z)$,
where $\alpha>1$, and put
$
y=(\alpha-1)/(\alpha+1).
$
Since
$
\tan\left(\frac{i}{2}\log\alpha\right)=iy,
$
the chain rule gives
\[
\left(\frac{i}{2}\right)^n T_n(iy)
=
\left(\alpha\frac{d}{d\alpha}\right)^n(iy).
\]
Moreover, since
$
y=1-2/(1+\alpha),
$
substituting $z=-\alpha$ into
\eqref{eq:ordinary-derivative} yields
\begin{equation}\label{eq:tangent-euler}
T_n(iy)
=-i^{1-n}2^{n+1}
\left(\alpha\frac{d}{d\alpha}\right)^n
\frac{1}{1+\alpha}
=i^{1-n}2^{n+1}
\frac{\alpha A_n(-\alpha)}
{(1+\alpha)^{n+1}}.
\end{equation}
Thus $A_n(-\alpha)=0$ implies $T_n(iy)=0$.

For the reciprocal root $-\alpha^{-1}$,
the corresponding value is
$
(\alpha^{-1}-1)/(\alpha^{-1}+1)=-y.
$
Applying the same identity with
$\alpha$ replaced by $\alpha^{-1}$
gives $T_n(-iy)=0$.

Next, solving
$
(1+t)/(1-t)=iy
$
gives
$
t_y:=(iy-1)/(iy+1).
$
Since $\alpha>1$, we have $0<y<1$.
Consequently,
\[
|t_y|=1,
\qquad
1+t_y^2
=\frac{2(1-y^2)}{(1+iy)^2}\ne0.
\]
As $t_y\ne1$, equation~\eqref{eq:mayeh}
implies that
\[
\Ahat_n(t_y)=0.
\]
The same argument applied to $-y$
shows that $\Ahat_n(t_{-y})=0$.

For each reciprocal pair, write
$\alpha=e^{v_i}$ with $v_i>0$.
The corresponding roots $t_y$ and $t_{-y}$
satisfy
$
t_{-y}=\overline{t_y},  	t_yt_{-y}=1.
$
Furthermore,
\[
t_y+t_{-y}
=\frac{2(y^2-1)}{1+y^2}
=-\frac{4\alpha}{\alpha^2+1}
=-2\sech v_i.
\]
Therefore, the corresponding quadratic factor is
\[
(t-t_y)(t-t_{-y})
=1+2\sech v_i\,t+t^2.
\]

When $n$ is even, the additional root $-1$
of $A_n$ corresponds to $y=0$.
Since $T_n(0)=0$ for even $n$,
equation~\eqref{eq:mayeh} gives
$\Ahat_n(-1)=0$.
This contributes the factor $1+t$.

Since $A_n$ has $n-1$ distinct negative roots
and the transformations above are injective,
we have obtained all $n-1$ roots of $\Ahat_n$.
Comparing constant terms gives~\eqref{eq:factor}.
Finally, $v_i>0$ implies
$0<c_i=2\sech v_i<2$,
completing the proof.
\end{proof}

A quadratic factor in~\eqref{eq:factor} is called
\emph{strong} if $c_i>3/2$.
Since $\sech v$ is decreasing for $v>0$,
we estimate the number of strong factors
by counting the zeros of $A_n(-e^v)$
in a suitable interval.

\begin{lemma}\label{lem:strong}
For $n\ge3$, at least
$\lfloor(n+1)/13\rfloor-1$
of the quadratic factors in~\eqref{eq:factor}
are strong.
\end{lemma}

\begin{proof}
We count zeros of $A_n(-e^v)$ in $0<v<7/9$ and then show that every such zero gives $c_i>3/2$.

\noindent\textbf{Step 1. Expressing the Eulerian polynomial as a sum over poles.}
Logarithmically differentiating the classical product
\[
\cosh\frac v2
=\prod_{j=0}^{\infty}
\left(1+\frac{v^2}{(2j+1)^2\pi^2}\right),
\]
and using the identity $(1+e^v)^{-1}=\tfrac12-\frac{d}{dv}\log\cosh(v/2)$ give
\begin{equation}\label{eq:poles}
\frac1{1+e^v}
=\frac12-\sum_{j=0}^{\infty}\left(
\frac1{v+(2j+1)\pi i}+
\frac1{v-(2j+1)\pi i}\right).
\end{equation}
The paired terms equal
\[
\frac{2v}{v^2+(2j+1)^2\pi^2}.
\]
The series and its derivatives of every positive
order converge locally uniformly on the real axis.
Thus, termwise differentiation is justified.

Define
\[
S_n(v)=\sum_{j\in\mathbb Z}
\frac{1}{(v+(2j+1)\pi i)^{n+1}}.
\]
Differentiating~\eqref{eq:poles} $n$ times gives
\[
\frac{d^n}{dv^n}\frac{1}{1+e^v}
=(-1)^{n+1}n!S_n(v).
\]
On the other hand, substituting $r=e^v$
into~\eqref{eq:ordinary-derivative} yields
\[
\frac{d^n}{dv^n}\frac{1}{1+e^v}
=-\frac{e^vA_n(-e^v)}{(1+e^v)^{n+1}}.
\]
Comparing these two identities, we obtain
\begin{equation}\label{eq:Sn}
A_n(-e^v)
=(-1)^n n!e^{-v}(1+e^v)^{n+1}S_n(v).
\end{equation}
Since the factor multiplying $S_n(v)$
is nonzero for real $v$, the two functions
$A_n(-e^v)$ and $S_n(v)$ have the same
real zeros.

\medskip
\noindent\textbf{Step 2. Estimating the contribution
of the remaining poles.}
For $0\le v\le7/9$, let
$
\theta(v)=\pi/2-\arctan (v/\pi).
$
The two terms in $S_n(v)$ corresponding
to $\pm\pi i$ give
\[
\frac1{(v+\pi i)^{n+1}}
+\frac1{(v-\pi i)^{n+1}}
=
\frac{2\cos((n+1)\theta(v))}
{(v^2+\pi^2)^{(n+1)/2}}.
\]
Consequently,
\[
\left|
\frac{(v^2+\pi^2)^{(n+1)/2}}2
S_n(v)
-\cos((n+1)\theta(v))
\right|
\le
\sum_{j=1}^{\infty}
\left(
\frac{\sqrt{v^2+\pi^2}}
{\sqrt{v^2+(2j+1)^2\pi^2}}
\right)^{n+1}.
\]
Since $0\le v\le7/9$ and $\pi>3$,
\[
\frac{\sqrt{v^2+\pi^2}}
{\sqrt{v^2+(2j+1)^2\pi^2}}
\le
\frac{\sqrt{1+(7/27)^2}}{2j+1}
<
\frac{21}{20(2j+1)}
<1.
\]
As $n+1\ge4$, we have
\begin{align*}
\sum_{j=1}^{\infty}
\left(
\frac{\sqrt{v^2+\pi^2}}
{\sqrt{v^2+(2j+1)^2\pi^2}}
\right)^{n+1}
&<
\left(\frac{21}{20}\right)^4
\sum_{j=1}^{\infty}\frac1{(2j+1)^4}\le
\left(\frac{21}{20}\right)^4
\left(
\frac1{81}
+\int_1^\infty\frac{dx}{(2x+1)^4}
\right)\\
&=
\left(\frac{21}{20}\right)^4
\frac1{54}
<\frac1{40}.
\end{align*}
Therefore,
\begin{equation}\label{eq:normalized-error}
\left|
\frac{(v^2+\pi^2)^{(n+1)/2}}2
S_n(v)
-\cos((n+1)\theta(v))
\right|
<\frac1{40}.
\end{equation}
In particular, whenever
$\cos((n+1)\theta(v))=\pm1$,
the function $S_n(v)$ has the same sign as the cosine.

\medskip
\noindent\textbf{Step 3. Counting the zeros by alternating signs.}
The function $\theta(v)$ is continuous
and strictly decreasing on $[0,7/9]$.
Hence $(n+1)\theta(v)/\pi$ runs through
the interval
\begin{equation}\label{eq:phase-interval}
\left[
\frac{n+1}{2}
-\frac{n+1}{\pi}\arctan\frac7{9\pi},
\frac{n+1}{2}
\right].
\end{equation}

We claim that this interval has length
greater than $(n+1)/13$.
Using $\pi<22/7$ and
$\arctan x>x-x^3/3$ for $x>0$, we obtain
\begin{equation}\label{eq:eta-bound}
\begin{aligned}
\frac1\pi\arctan\frac7{9\pi}
>\frac7{22}\arctan\frac{49}{198}
>\frac7{22}
\left(
\frac{49}{198}
-\frac13\left(\frac{49}{198}\right)^3
\right)
>\frac1{13}.
\end{aligned}
\end{equation}

Thus the interval~\eqref{eq:phase-interval}
contains at least
$\lfloor(n+1)/13\rfloor$
consecutive integers.
For each such integer $k$, there is a unique
$v\in[0,7/9]$ satisfying
$
(n+1)\theta(v)=k\pi.
$
At these points,
\[
\cos((n+1)\theta(v))=(-1)^k.
\]
By~\eqref{eq:normalized-error},
the corresponding values of $S_n(v)$
are nonzero and alternate in sign.

The intermediate value theorem therefore
gives a zero of $S_n(v)$ between each
pair of consecutive points.
These zeros are distinct and lie in $(0,7/9)$.
Hence $S_n(v)$ has at least
\[
\max\left\{
0,\left\lfloor\frac{n+1}{13}\right\rfloor-1
\right\}
\]
distinct zeros in $(0,7/9)$.
By~\eqref{eq:Sn}, the same holds
for $A_n(-e^v)$.

\medskip
\noindent\textbf{Step 4. Obtaining strong factors.}
It remains to show that every zero counted
above corresponds to a strong factor.
Using the Taylor series for $\cosh v$
and bounding its terms from $v^2/2$
onward by a geometric series, we obtain
\begin{equation}\label{eq:cosh-bound}
\cosh\frac79
\le
1+\frac{(7/9)^2/2}
{1-(7/9)^2/12}
=\frac{1217}{923}
<\frac43.
\end{equation}
Thus
$
2\sech(7/9)>3/2.
$
Since $\sech v$ is strictly decreasing
for $v>0$, every zero $v_i\in(0,7/9)$
satisfies
\[
c_i=2\sech v_i>\frac32.
\]
By Proposition~\ref{prop:factor},
distinct positive zeros $v_i$ correspond
to distinct quadratic factors.
Therefore at least
$\lfloor(n+1)/13\rfloor-1$
factors are strong, as required. 
This proves the lemma.
\end{proof}

\section{Log-concavity under parameter variation}\label{sec:tools}
To apply the strong-factor estimate from Section~2, we establish two tools for products of reciprocal quadratic factors.

Let $q(t)=\sum_{j=0}^{d}b_jt^j$ with \(b_j=0\) for \(j<0\) or \(j>d\).
For each $k\in\mathbb Z$, the identity
\[
[t^k]\bigl(q(t)(1+xt+t^2)\bigr)
=b_k+xb_{k-1}+b_{k-2}
\]
gives
\begin{equation}\label{eq:quadratic-delta}
\Delta_k\bigl(q(t)(1+xt+t^2)\bigr)
=\Delta_k\bigl(q(t)(1+t^2)\bigr)+L_k(q)x+\Delta_{k-1}(q)x^2,
\end{equation}
where
$
L_k(q)
=(b_kb_{k-1}-b_{k+1}b_{k-2})
+(b_{k-1}b_{k-2}-b_kb_{k-3}).
$
By Lemma~\ref{lem:mixed}, $L_k(q)\ge0$
whenever $q(t)$ is a product of factors
$1+c_it+t^2$ with $c_i\ge0$. The preceding identity leads to the following
endpoint reduction criterion for log-concavity.

\begin{proposition}\label{prop:endpoint}
Let $m\ge0$ and $0\le a_i\le d_i$ for
$1\le i\le m$. For each fixed integer $k$,
the minimum of
\[
\Delta_k\left(
\prod_{i=1}^{m}(1+c_it+t^2)
\right)
\]
over $a_i\le c_i\le d_i$ is attained when $c_i\in\{a_i,d_i\}$ for all $i$.
Consequently, the coefficient sequence of the product is log-concave for all	$a_i\le c_i\le d_i$ if and only if	it is log-concave whenever $c_i\in\{a_i,d_i\}$ for all $i$.
\end{proposition}

\begin{proof}
The case $m=0$ is immediate.
Fix $k$ and all parameters except $c_i=x$,
and let
\[
q(t)=\prod_{j\ne i}(1+c_jt+t^2).
\]
By~\eqref{eq:quadratic-delta}, we have
$
\Delta_k\bigl(q(t)(1+xt+t^2)\bigr)
=\Delta_k\bigl(q(t)(1+t^2)\bigr)
+L_k(q)x+\Delta_{k-1}(q)x^2.
$
Since $L_k(q)\ge0$ by
Lemma~\ref{lem:mixed}, we consider
two cases.

If $\Delta_{k-1}(q)\ge0$, then
\[
\frac{d}{dx}
\Delta_k\bigl(q(t)(1+xt+t^2)\bigr)
=L_k(q)+2x\Delta_{k-1}(q)\ge0
\]
for $x\ge0$.
Thus the minimum on $[a_i,d_i]$
is attained at $a_i$.

If $\Delta_{k-1}(q)<0$, then
\eqref{eq:quadratic-delta} is a concave quadratic in $x$.
Its minimum on $[a_i,d_i]$ is therefore attained
at an endpoint.

Thus, for fixed values of the other parameters,
the minimum with respect to $c_i$ is attained
at one of the endpoints of $[a_i,d_i]$.

Since the expression is continuous on the
compact parameter box, it attains a global
minimum. Starting from a minimizer, we can
successively replace each parameter by a
suitable endpoint without increasing the value.
After all parameters have been replaced,
we obtain a vertex at which the global
minimum is attained.

Consequently, if the coefficient sequence
of the product is log-concave at every vertex, then $\Delta_k\ge0$ throughout the box for every relevant $k$.
The converse is immediate.
\end{proof}

Suppose that the coefficient sequence of $q(t)$
is nonnegative, log-concave, and has no internal
zeros. Then
\begin{equation}\label{eq:mixed-general}
b_kb_{k-1}\ge b_{k+1}b_{k-2}
\qquad(k\in\mathbb Z).
\end{equation}
If $b_{k+1}b_{k-2}=0$, the inequality is
immediate. Otherwise, the absence of internal zeros
implies that $b_{k-1}b_k>0$.
Multiplying the two log-concavity inequalities
\[
b_{k-1}^2\ge b_{k-2}b_k,
\qquad
b_k^2\ge b_{k-1}b_{k+1},
\]
and dividing by $b_{k-1}b_k$ gives
\eqref{eq:mixed-general}.
Applying this inequality at $k$ and $k-1$
gives $L_k(q)\ge0$, as defined
in~\eqref{eq:quadratic-delta}.

The following lemma allows us to replace factors $1+t^2$ by $1+xt+t^2$ while preserving log-concavity under suitable assumptions.

\begin{lemma}\label{lem:smoothing}
Let $g(t)$ be a nonzero polynomial with
nonnegative coefficients and no internal zeros,
and let $N$ be a nonnegative integer.
Suppose that the coefficient sequence of
\[
g(t)(1+t^2)^j
\]
is log-concave for every $0\le j\le N$.
Then, for every integer \(0\le s\le N\) and
$x_1,\ldots,x_s\ge0$, the coefficient sequence of
\[
g(t)\prod_{i=1}^{s}(1+x_it+t^2)
\]
is log-concave.
\end{lemma}

\begin{proof}
First, suppose that $q(t)$ has nonnegative,
log-concave coefficients and no internal zeros,
and that $q(t)(1+t^2)$ also has a log-concave
coefficient sequence.
By~\eqref{eq:quadratic-delta},
\[
\Delta_k\bigl(q(t)(1+xt+t^2)\bigr)=
\Delta_k\bigl(q(t)(1+t^2)\bigr)
+L_k(q)x+\Delta_{k-1}(q)x^2.
\]
All three terms on the right are nonnegative
for $x\ge0$, since $L_k(q)\ge0$
by~\eqref{eq:mixed-general}.
Consequently,
\begin{equation}\label{eq:one-factor}
q(t)(1+xt+t^2)
\end{equation}
has a log-concave coefficient sequence
for every $x\ge0$.

If $N=0$, the conclusion follows directly from the hypothesis.
We assume for the remainder of the proof that $N\ge1$.
Since $g(t)(1+t^2)$ has a log-concave	coefficient sequence, $g(t)$ cannot be a monomial. Indeed, if $g(t)=at^r$
with $a>0$, then
$
\Delta_{r+1}\bigl(g(t)(1+t^2)\bigr)
=-a^2<0,
$
a contradiction.
Thus the support of $g(t)$ contains at least two consecutive indices.
Multiplication by $1+xt+t^2$, where $x\ge0$, preserves this property and introduces no internal zeros.
Hence all products considered below have nonnegative coefficients and no internal zeros.

We prove by induction on $\ell$ that
\begin{equation}\label{eq:smoothing-induction}
g(t)(1+t^2)^j
\prod_{i=1}^{\ell}(1+x_it+t^2)
\end{equation}
has a log-concave coefficient sequence
whenever $j,\ell\ge0$ and $j+\ell\le N$.

For $\ell=0$, the assertion follows
from the hypothesis of the lemma.
Suppose that it holds for $\ell-1$,
and fix $j\ge0$ with $j+\ell\le N$.

Since
$
j+\ell-1\le N,
(j+1)+(\ell-1)=j+\ell\le N,
$
the induction hypothesis shows that both
\[
g(t)(1+t^2)^j
\prod_{i=1}^{\ell-1}(1+x_it+t^2)
\quad
\text{ and }
\quad
g(t)(1+t^2)^{j+1}
\prod_{i=1}^{\ell-1}(1+x_it+t^2)
\]
have log-concave coefficient sequences.
Moreover, the first product has nonnegative
coefficients and no internal zeros,
as established above.
Since the second product is the first
multiplied by $1+t^2$,
applying~\eqref{eq:one-factor}
with $x=x_\ell$ shows that
\[
g(t)(1+t^2)^j
\prod_{i=1}^{\ell}(1+x_it+t^2)
\]
has a log-concave coefficient sequence.
This completes the induction.
Taking $j=0$ and $\ell=s$
gives the desired conclusion.
\end{proof}

\section{Proof of Theorem~\ref*{thm:main}}

We first establish a block criterion for log-concavity
and then combine it with Lemma~\ref{lem:strong}
to prove Theorem~\ref{thm:main}.

\begin{proposition}\label{prop:finite-block}
For integers $0\le s\le73$ and $0\le j\le8$,
let
\begin{equation}\label{eq:P-sj}
P_{s,j}(t)
=(1+t^2)^s(2+3t+2t^2)^{8-j}
(2+4t+2t^2)^j.
\end{equation}
Then
\begin{equation}\label{eq:finite-block}
\Delta_k(P_{s,j})\ge0,
\qquad 1\le k\le8+s.
\end{equation}
Equality holds only when $(s,j,k)=(73,0,1)$.
\end{proposition}

\begin{proof}
Each factor in~\eqref{eq:P-sj} is palindromic, so
$P_{s,j}(t)$ is palindromic of degree $2(8+s)$ and
$
\Delta_k(P_{s,j})=\Delta_{2(8+s)-k}(P_{s,j}).
$
Thus the inequalities in~\eqref{eq:finite-block}
imply log-concavity of the entire coefficient sequence.
We verify these inequalities as follows.

For each $0\le j\le8$, we first compute the coefficients of $P_{0,j}(t)$ by multiplying
its eight quadratic factors. The coefficients of $P_{s+1,j}(t)$ are then obtained recursively from
\begin{equation}\label{eq:block-rec}
[t^k]P_{s+1,j}(t)
=[t^k]P_{s,j}(t)+[t^{k-2}]P_{s,j}(t),
\qquad 0\le s<73.
\end{equation}
Using these coefficients, we evaluate
$\Delta_k(P_{s,j})$ for
$1\le k\le8+s$.

The verification involves
$
9\sum_{s=0}^{73}(8+s)=29\,637
$
integer comparisons. All the computed values are nonnegative, with equality occurring only at $(s,j,k)=(73,0,1)$.

The accompanying script
\texttt{verify\_log\_concavity.py},
available in the
\href{https://github.com/lijiewang86/alternating-eulerian-log-concavity}
{GitHub repository},
reproduces this verification using exact integer arithmetic.
This completes the proof.
\end{proof}

We now combine Proposition~\ref{prop:finite-block}
with the results of Section~3 to obtain
the following block criterion.

\begin{proposition}\label{thm:block}
Suppose $c_1,\dots,c_8\in [3/2,2]$, $s\in\mathbb Z$ with $0\le s\le 73$, and $x_1,\dots,x_s\ge 0$.
Then the coefficient sequence of the polynomial
\begin{equation}\label{eq:block-product}
\prod_{i=1}^{8}(1+c_it+t^2)
\prod_{\ell=1}^{s}(1+x_\ell t+t^2)
\end{equation}
is log-concave.
\end{proposition}

\begin{proof}
Let
$
g(t)=\prod_{i=1}^{8}(1+c_it+t^2).
$
We first show that $g(t)(1+t^2)^s$
has a log-concave coefficient sequence
for every integer $0\le s\le73$.
Fix such an $s$.
By Proposition~\ref{prop:endpoint},
it suffices to consider
$c_i\in\{3/2,2\}$ for $1\le i\le8$,
with the remaining $s$ parameters
fixed at zero.
If exactly $j$ of the eight parameters
equal $2$, the corresponding product is
\[
g(t)(1+t^2)^s
=(1+t^2)^s
\left(1+\frac32t+t^2\right)^{8-j}
(1+2t+t^2)^j=2^{-8}P_{s,j}(t).
\]
By Proposition~\ref{prop:finite-block}
and the palindromicity of $P_{s,j}(t)$,
its coefficient sequence is log-concave
for every $0\le j\le8$.
Since multiplication by a positive
constant preserves log-concavity,
the same holds for $2^{-8}P_{s,j}(t)$.
Proposition~\ref{prop:endpoint} therefore
shows that $g(t)(1+t^2)^s$ has a
log-concave coefficient sequence
for all $c_i\in[3/2,2]$.

Moreover, $g(t)$ has positive coefficients
and hence has no internal zeros.
Thus all the hypotheses of
Lemma~\ref{lem:smoothing} are satisfied
with $N=73$.
It follows that
\[
g(t)\prod_{\ell=1}^{s}(1+x_\ell t+t^2)
\]
has a log-concave coefficient sequence
for every $0\le s\le73$ and
$x_1,\ldots,x_s\ge0$,
as required.
\end{proof}

\begin{remark}\label{rem:sharp}
The bound $s\le73$ in Proposition~\ref{thm:block}
is sharp. Taking $c_i=3/2$ for
$1\le i\le8$ and $x_\ell=0$ for
$1\le\ell\le s$ gives
\[
\left(1+\frac32t+t^2\right)^8(1+t^2)^s.
\]
Its first three coefficients are $1$, $12$, and
$
8+\binom82\left(\frac32\right)^2+s=71+s,
$
respectively. Hence the log-concavity
difference at index $1$ is
$
12^2-(71+s)=73-s,
$
which is negative when $s=74$.
Thus the bound cannot be improved.
\end{remark}

We now complete the proof of
Theorem~\ref{thm:main} using exact integer verification for $n\le220$ and the strong-factor estimate together with the block criterion for $n\ge221$.

\begin{proof}[\textbf{Proof of Theorem~\ref{thm:main}}]
\noindent\textbf{Case 1. $1\le n\le220$.}
We use the recurrence of Lin, Ma, Wang and
Wang~\cite[Theorem~2.1]{LinMaWangWang2022}:
\begin{equation}\label{eq:alt-rec}
2a_{n+1,k}
=(k+1)(a_{n,k+1}+a_{n,k-1})+(n-k+1)(a_{n,k}+a_{n,k-2}),
\qquad a_{1,0}=1,
\end{equation}
where \(a_{n,k}=0\) unless \(0\le k\le n-1\). 
Starting from $a_{1,0}=1$, this recurrence
determines the coefficients of $\Ahat_n(t)$
for every $n\le220$.

By~\eqref{eq:symmetry}, we have
$
\Delta_k(\Ahat_n)
=\Delta_{n-1-k}(\Ahat_n).
$
Thus it suffices to verify
$\Delta_k(\Ahat_n)\ge0$ for
$1\le k\le\lfloor(n-1)/2\rfloor$.
The number of required comparisons is
$
\sum_{n=1}^{220}
\left\lfloor\frac{n-1}{2}\right\rfloor
=11\,990.
$
An exact integer computation using
\eqref{eq:alt-rec} shows that all these
differences are nonnegative, with equality
only at $(n,k)=(3,1)$.
The accompanying script
\texttt{verify\_log\_concavity.py}
reproduces this verification.
Hence the theorem holds for $n\le220$.

\medskip
\noindent\textbf{Case 2. $n\ge221$.}
Let $m=\lfloor(n-1)/2\rfloor$.
By Proposition~\ref{prop:factor}, we have
\[
\Ahat_n(t)
=a_{n,0}(1+t)^{n-1-2m}
\prod_{i=1}^{m}(1+c_it+t^2),
\]
where $a_{n,0}>0$ and $0<c_i<2$.

Let
$
b=\left\lfloor(2m-11)/104\right\rfloor.
$
Since $m\ge110$, we have $b\ge2$.
By Lemma~\ref{lem:strong}, the number
of strong factors is at least
\begin{equation}\label{eq:strong-main}
\left\lfloor\frac{n+1}{13}\right\rfloor-1
\ge
\left\lfloor\frac{2m-11}{13}\right\rfloor\ge8b.
\end{equation}
Moreover, the definition of $b$ gives
$
2m-11<104(b+1),
$
and hence
\begin{equation}\label{eq:block-number}
m\le52b+57\le81b,
\end{equation}
where the last inequality follows from
$b\ge2$.

We can therefore choose $8b$ strong
factors and divide them into $b$ groups,
each containing eight strong factors.
By~\eqref{eq:block-number}, the number
of remaining factors satisfies
\[
m-8b\le73b.
\]
These factors can be distributed among
the $b$ groups so that each group contains
at most $73$ additional factors.

The eight strong factors in each group
have parameters in $[3/2,2]$, while
all remaining parameters are nonnegative.
Thus, by Proposition~\ref{thm:block},
every group has a log-concave coefficient
sequence. Since all the factors have
positive coefficients, each group also
has no internal zeros.

Applying Lemma~\ref{lem:convolution}
repeatedly, we conclude that
\[
\prod_{i=1}^{m}(1+c_it+t^2)
\]
has a log-concave coefficient sequence
with no internal zeros.

Finally, $n-1-2m\in\{0,1\}$.
Since the coefficient sequence of $1+t$
is log-concave, another application of
Lemma~\ref{lem:convolution}
shows that multiplying by the possible
factor $1+t$ preserves log-concavity.
Multiplication by $a_{n,0}>0$ also
preserves log-concavity.
Therefore, the coefficient sequence of
$\Ahat_n(t)$ is log-concave for all
$n\ge221$, which completes the proof.
\end{proof}
As a consequence of Theorem~\ref{thm:main}, we recover
the unimodality established by Lin, Ma, Wang and
Wang~\cite{LinMaWangWang2022}.

\begin{corollary}\label{cor:unimodality}
For every $n\ge1$, the coefficient sequence of
$\Ahat_n(t)$ is unimodal.
\end{corollary}

\section{Applications}
\label{sec:consequences}

We give four consequences of Theorem~\ref{thm:main}.
We first transfer log-concavity to the distribution of
$3$-descents and derive a local limit theorem.
We then bound the largest coefficient and obtain a
Lorentzian reformulation of the main theorem.


For $\pi\in\mathfrak S_{n+1}$, let $d_3(\pi)$ denote the number of indices $i\in[n-1]$ for which the consecutive triple $\pi_i\pi_{i+1}\pi_{i+2}$ has relative order $132$, $213$ or $321$.
Chebikin's equidistribution theorem
\cite[Theorem~2.3 and Corollary~2.2]{Chebikin2008}
relates this statistic to alternating descents.

\begin{corollary}\label{cor:three-descents}
For every $n\ge1$,
\begin{equation}\label{eq:three-descent-enumerator}
\sum_{\pi\in\mathfrak S_{n+1}}t^{d_3(\pi)}
=(n+1)\Ahat_n(t).
\end{equation}
Consequently, the coefficient sequence on the left is log-concave for every $n\ge1$.
\end{corollary}


Normal approximation for alternating descents was studied
by Heerten, Sambale and
Th\"ale~\cite{HeertenSambaleThale2024}.
Theorem~\ref{thm:main} supplies the log-concavity hypothesis
in Bender's local limit theorem~\cite{Bender1973}.
We give the resulting coefficient estimate below.

\begin{corollary}\label{cor:local-limit}
As $n\to\infty$,
\begin{equation}\label{eq:local-limit}
\sup_{k\in\mathbb Z}
\left|
\sqrt{\frac{5n}{12}}\frac{a_{n,k}}{n!}
-\frac{1}{\sqrt{2\pi}}
\exp\left(
-\frac{6}{5n}\left(k-\frac{n-1}{2}\right)^2
\right)
\right|\longrightarrow0,
\end{equation}
where $a_{n,k}=0$ outside $0\le k\le n-1$.
Consequently,
\begin{equation}\label{eq:central-asymptotic}
a_{n,\lfloor(n-1)/2\rfloor}
\sim n!\sqrt{\frac{6}{5\pi n}}.
\end{equation}
\end{corollary}
\begin{proof}
Hwang, Chern and Duh
\cite[Section~5.4.2]{HwangChernDuh2020}
established asymptotic normality with centering $n/2$
and scale $\sqrt{5n/12}$.
The centering may be replaced by $(n-1)/2$, since
the difference divided by $\sqrt{5n/12}$ tends to zero.
By Proposition~\ref{prop:factor}, the coefficients are positive,
and by Theorem~\ref{thm:main}, they form a log-concave sequence.
Since $\sqrt{5n/12}\to\infty$, Bender's local limit
theorem~\cite{Bender1973} gives~\eqref{eq:local-limit}.

Taking $k=\lfloor(n-1)/2\rfloor$ in~\eqref{eq:local-limit},
the exponential term tends to $1$, giving
\eqref{eq:central-asymptotic}.
\end{proof}

The main theorem also yields bounds for the largest coefficient through a concentration inequality for symmetric log-concave distributions.

\begin{corollary}\label{cor:coefficient-bounds}
For every $n\ge2$, the maximum coefficient of $\Ahat_n(t)$ is
$a_{n,\lfloor(n-1)/2\rfloor}$, and
\begin{equation}\label{eq:coefficient-bounds}
\frac{n!}{\sqrt{5n-6}}
\le a_{n,\lfloor(n-1)/2\rfloor}
\le n!\sqrt{\frac{6}{5n-1}}.
\end{equation}
\end{corollary}
\begin{proof}
Palindromicity and Corollary~\ref{cor:unimodality} show that
the central coefficient is a maximum coefficient.
Let $X_n$ be the number of alternating descents in a uniformly
random permutation in $\mathfrak S_n$.
Then $\mathbb P(X_n=k)=a_{n,k}/n!$.

Each alternating descent indicator has variance $1/4$.
Two adjacent indicators are simultaneously $1$ with probability
$1/3$, so their covariance is $1/12$; nonadjacent indicators
are independent. Hence
\[
\operatorname{Var}(X_n)
=\frac{n-1}{4}+\frac{n-2}{6}
=\frac{5n-7}{12}.
\]
By Theorem~\ref{thm:main} and Proposition~\ref{prop:factor},
the distribution of $X_n$ is log-concave with no internal zeros,
and palindromicity makes it symmetric about $(n-1)/2$.
Thus Bobkov, Marsiglietti and Melbourne
\cite[Theorem~1.1]{BobkovMarsigliettiMelbourne2022} give
\[
\frac{1}{\sqrt{1+12\operatorname{Var}(X_n)}}
\le \max_k\frac{a_{n,k}}{n!}
\le \frac{1}{\sqrt{1+2\operatorname{Var}(X_n)}}.
\]
Substituting the variance above and multiplying by $n!$
gives~\eqref{eq:coefficient-bounds}.
\end{proof}

We finish with a reformulation of Theorem~\ref{thm:main}
using the characterization of bivariate Lorentzian polynomials.

\begin{corollary}\label{cor:lorentzian}
For every $n\ge1$, the homogeneous polynomial
\begin{equation}\label{eq:lorentzian-polynomial}
\sum_{k=0}^{n-1}\binom{n-1}{k}a_{n,k}x^ky^{n-1-k}
\end{equation}
is Lorentzian and completely log-concave.
\end{corollary}

\begin{proof}
The cases $n=1,2$ are immediate.
For $n\ge3$, the bivariate characterization in
\cite[Example~2.26]{BrandenHuh2020} reduces the Lorentzian
property to log-concavity and the absence of internal zeros
of the coefficients after division by $\binom{n-1}{k}$.
For~\eqref{eq:lorentzian-polynomial}, these normalized
coefficients are exactly $a_{n,k}$.
They are positive by Proposition~\ref{prop:factor}
and log-concave by Theorem~\ref{thm:main}.
Complete log-concavity follows from
\cite[Theorem~2.30]{BrandenHuh2020}.
\end{proof}

\section*{Conflict of Interest}  
The authors declare no conflict of interest.

\section*{Acknowledgements}
This work was supported in part by the National Natural Science Foundation of China (Grant Nos.~12271078, 12101110) and the Scientific Research Fund of Liaoning Provincial Education Department (Grant No.~LJ212410173071).

\end{document}